%

%

\def\today{\ifcase\month\or January\or February\or March\or
April\or May\or June\or July\or August\or September\or
October\or November\or December\fi \space\number\day,
\number\year}



\def\dspace{\lineskip=2pt\baselineskip=18pt\lineskiplimit=0pt}

\font \bbrm=cmbx10 at 12pt

\def\bigtype{\bbrm}

\hsize=13.5cm
\magnification=1200
\def\ce{\centerline}

\def\hb{\hfill\break}

\def\title #1{\null\bigskip\ce{\bigtype #1}\bigskip}

\def\alp{\alpha}		
\def\bet{\beta}		
\def\gam{\gamma}		
\def\del{\delta}

\def\kap{\kappa}
\def\lam{\lambda}		
\def\sig{\sigma}		

\def\ome{\omega}		


\def\calF{{\cal F}}
\def\calG{{\cal G}}

\def\calK{{\cal K}}
\def\calL{{\cal L}}

\def\calP{{\cal P}}

\def\calU{{\cal U}}



    
\font\tenboldgreek=cmmib10  \font\sevenboldgreek=cmmib10 at
7pt
\font\fiveboldgreek=cmmib10 at 7pt
\newfam\bgfam
\textfont\bgfam=\tenboldgreek \scriptfont\bgfam=\sevenboldgreek
\scriptscriptfont\bgfam=\fiveboldgreek

\mathchardef\ggarrow="7010

\font\tengerman=eufm10 \font\sevengerman=eufm7 \font\fivegerman=eufm5
\font\tendouble=msym10 \font\sevendouble=msym7 \font\fivedouble=msym5

\textfont4=\tengerman \scriptfont4=\sevengerman
\scriptscriptfont4=\fivegerman
\newfam\dbfam
\textfont\dbfam=\tendouble \scriptfont\dbfam=\sevendouble
\scriptscriptfont\dbfam=\fivedouble
\def\gr{\fam4}

\mathchardef\ng="702D
\mathchardef\dbA="7041
\mathchardef\sm="7072
\mathchardef\nvdash="7030
\mathchardef\nldash="7031
\mathchardef\lne="7008
\mathchardef\sneq="7024
\mathchardef\spneq="7025
\mathchardef\sne="7028
\mathchardef\spne="7029
\mathchardef\ltms="706E
\mathchardef\tmsl="706F

\mathchardef\dbA="7041

\def\supsetneqq{\,{\fam=\dbfam\spneq}\,}

	\def\gra{{\gr a}}

\mathchardef\dbA="7041 
\mathchardef\dbB="7042 
\mathchardef\dbC="7043 
\mathchardef\dbD="7044 
\mathchardef\dbE="7045 
\mathchardef\dbF="7046 
\mathchardef\dbG="7047 
\mathchardef\dbH="7048 
\mathchardef\dbI="7049 
\mathchardef\dbJ="704A 
\mathchardef\dbK="704B 
\mathchardef\dbL="704C 
\mathchardef\dbM="704D 
\mathchardef\dbN="704E 
\mathchardef\dbO="704F 
\mathchardef\dbP="7050 
\mathchardef\dbQ="7051 
\mathchardef\dbR="7052 
\mathchardef\dbS="7053 
\mathchardef\dbT="7054 
\mathchardef\dbU="7055 \def\UU{{\fam=\dbfam\dbU}}
\mathchardef\dbV="7056 
\mathchardef\dbW="7057 
\mathchardef\dbX="7058 
\mathchardef\dbY="7059 
\mathchardef\dbZ="705A 

\def\nek{,\ldots,}
\def\sdp{\times \hskip -0.3em {\raise 0.3ex
\hbox{$\scriptscriptstyle |$}}} 


\def\dom{\mathop{\rm dom}\nolimits}

\def\min{\mathop{\rm min}}

\def\supp{\mathop{\rm supp}}


\def\oB{{\overline B}}

\def\oC{{\overline C}}

\def\oh{{\overline h}}

\def\om{{\overline m}}


\def\obet{{\overline\bet}}

\def\odel{{\overline\delta}}

\def\ogam{{\overline\gam}}

\def\okap{{\overline\kap}}

\def\onu{{\overline\nu}}

\def\orho{{\overline\rho}}

\def\otau{{\overline\tau}}

\def\oxi{{\overline\xi}}







\def\sqr#1#2{{\vcenter{\hrule height.#2pt\hbox{\vrule
width.#2pt height#1pt \kern#1pt \vrule width.#2pt}\hrule
height.#2pt}}}
\def\square{\mathchoice{\sqr34}{\sqr34}{\sqr{2.1}3}{\sqr{1.5}3}}

\def\buildrul#1\under#2{\mathrel{\mathop{\null#2}\limits_{#1}}}

\def\boxit#1{\vbox{\hrule\hbox{\vrule\kern3pt\vbox{\kern3pt#1
\kern3pt}\kern3pt\vrule}\hrule}}

\def\prodl{\prod\limits}

\def\bigcupl{\bigcup\limits}

\def\subheading#1{\medskip\goodbreak\noindent{\bf
#1.}\quad}
\def\sect#1{\goodbreak\bigskip\centerline{\bf#1}\medskip}
\def\pr{\smallskip\noindent{\bf Proof:\quad}}
\def\onumber #1{\ooalign{\hfil\raise.07ex\hbox{\hfill$\scriptstyle
\,#1$\hfil}
\cr\cr{$\bigcirc$}}}
\def\onumber c{\ooalign{\hfil\raise.07ex\hbox{\hfill$\scriptstyle
\,c$\hfil}
\cr\cr{$\bigcirc$}}}
\def\alpcirc {\ooalign{\hfil\raise.07ex
\hbox{\hfill$\scriptstyle\alp\;$\hfill}\cr\cr{$\bigcirc$}}}

\def\longmapright #1 #2 {\smash{\mathop{\hbox to
#1pt {\rightarrowfill}}\limits^{#2}}}
\def\longmapleft #1 #2 {\smash{\mathop{\hbox to
#1pt {\leftarrowfill}}\limits^{#2}}}

\def\references#1{\goodbreak\bigskip\par\centerline{\bf
References}\medskip\parindent=#1pt}
\def\ref#1{\par\smallskip\hang\indent\llap{\hbox
to \parindent{#1\hfil\enspace}}\ignorespaces}

\def\back{{\raise 2.5pt\hbox{$\,\scriptscriptstyle\backslash\,$}}}
\def\bks{{\backslash}}
\def\part{\partial}
\def\lwr #1{\lower 5pt\hbox{$#1$}\hskip -3pt}
\def\rse #1{\hskip -3pt\raise 5pt\hbox{$#1$}}
\def\lwrs #1{\lower 4pt\hbox{$\scriptstyle #1$}\hskip -2pt}
\def\rses #1{\hskip -2pt\raise 3pt\hbox{$\scriptstyle #1$}}

\def\<#1{\left\langle{#1}\right\rangle}

\def\subinbn{{\subset\hskip-8pt\raise 0.95pt\hbox{$\scriptscriptstyle\subset$}}}

\def\llvdash{\mathop{\|\hskip-2pt \raise 3pt\hbox{\vrule
height 0.25pt width 1.5cm}}}

\def\lvdash{\mathop{|\hskip-2pt \raise 3pt\hbox{\vrule
height 0.25pt width 1.5cm}}}

\def\fakebold#1{\leavevmode\setbox0=\hbox{#1}%
  \kern-.025em\copy0 \kern-\wd0
  \kern .025em\copy0 \kern-\wd0
  \kern-.025em\raise.0333em\box0 }

\font\msxmten=msxm10
\font\msxmseven=msxm7
\font\msxmfive=msxm5
\newfam\myfam
\textfont\myfam=\msxmten
\scriptfont\myfam=\msxmseven
\scriptscriptfont\myfam=\msxmfive
\mathchardef\rhookupone="7016
\mathchardef\ldh="700D
\mathchardef\leg="7053
\mathchardef\ANG="705E
\mathchardef\lcu="7070
\mathchardef\rcu="7071
\mathchardef\leseq="7035
\mathchardef\qeeg="703D
\mathchardef\qeel="7036
\mathchardef\blackbox="7004
\mathchardef\bbx="7003
\mathchardef\simsucc="7025

\def\rhookup{{\fam=\myfam \rhookupone}}

\font\tencaps=cmcsc10
\def\smallcaps{\tencaps}

\def\author#1{\bigskip\ce{\smallcaps #1}\medskip}

\def\tagg{^{\prime\prime}}
\def\taggg{^{\prime\prime\prime}}

\def\upddots{\mathinner{\mkern
1mu\raise 1pt \hbox{.}\mkern 2mu \mkern
2mu \raise 4pt\hbox{.}\mkern 1mu \raise 7pt\vbox {\kern 7
pt\hbox{.}}} }

\overfullrule=0pt
\def\lft{{\longleftrightarrow}}
{\nopagenumbers
\null
\vskip 4 true cm
\title{ON HIDDEN EXTENDERS}
\sect{Moti Gitik} \vskip 1 true cm
\ce{School of Mathematics}
\ce{Raymond and Beverly Sackler}
\ce{Faculty of Exact Sciences}
\ce{Tel Aviv University}
\ce{Tel Aviv 69978 Israel}
\vfill\eject}
\count0=1
\null
\dspace
\sect{1.~~Introduction}

Let $\kap$  be a singular cardinal violating
GCH  or a measurable with $2^\kap >\kap^+$.
The strength of this hypotheses was studied in
[Git1,2] and [Git-Mit] combining Shelah's $pcf$
theory with Mitchell's covering lemma.  The
basic approach was to use the covering lemma in
order to change sequences witnessing $pcf$ to
better and better one.    

A principal point there was to reconstruct
extenders as in [Git-] or to show that extenders
of an indiscernible sequence do not depend on a
particular precovering set as in [Git-Mit].  This
was shown to be true for $cof\kap >\ome$  or under
the assumption that for some $n<\ome$  $\ \{\alp <\kap |
\calK\models o(\alp)\ge\alp^{+n}\}$  is bounded in
$\kap$.	

The purpose of this paper will be to show the above
breaks down if for every $n<\ome\ \{\alp
<\kap|\calK \models o(\alp)\ge\alp^{+n}\}$ 
is cofinal in $\kap$.

\sect{2.~~Types of Ordinals}

Fix $n<\ome$.  For $k\le n$  let us consider a
language $\calL_{n,k}$  containing a constant
$c_\alp$ for every $\alp<\kap_n^{+k}$  and a
structure $\gra_{n,k}=\langle
H(\lam^{+k)}),\in,\lam,\le,0,1\ldots
,\alp,\ldots |\alp <\kap_n^{+k}\rangle$ in this
language, where $\lam$ is a regular cardinal big
enough.  For an ordinal $\xi <\lam$  (usually
$\xi <\kap_n^{+n+2}$) denote by $tp_{n,k}(\xi)$
the $\calL_k$-type realized by $\xi$  in
$\gra_{n,k}$.  Further we shall drop $n$
whenever it will not lead to confusion.  If $\del$
is an ordinal below $\lam$  then we shall consider a
language $\calL_{k,\del}$ obtained from
$\calL_k$  by adding a new constant $c$  and a
structure $\gra_{k,\del}$  as above with $c$
interpreted as $\del$.  $tp_{k,\del}(\xi)$ is
defined in the obvious fashion.      

\proclaim Lemma 2.1.  Suppose $2<k\le n$  and
$tp_k(\gam)=tp_k(\del)$.  Then for every $\xi
<\lam$  there is $\rho <\lam$  such that
$tp_{k-1,\gam}(\xi)=tp_{k-1,\del}(\rho)$.

\pr $tp_{k-1,\gam}(\xi)$  can be viewed as
a subset of $\kap_n^{+k}$  and by GCH as an
ordinal less than $\kap_n^{+k}$. 
Since $H(\lam^{+(k-1)}) \in H(\lam^{+k})$, 
the existence of $\xi$ having $tp_{k-1,\gam}(\xi)$ 
can be expressed by single true in $\gra_k$  formula
of $\calL_k$.  Hence $\del$  satisfies this formula
as well.  So there is $\rho <\lam$ with
$tp_{k-1,\del}(\rho)=tp_{k-1,\gam}(\xi)$.\hfill$\square$ 

The restriction $k>2$  is not essential for the
proof but further this will be the only
interesting case.  Since we are going to deal
with measures over $\kap_n$'s, the total number
of them is $\kap^{+2}_n$.  So $k>2$  will be
responsible for not dropping to isolated types. 

It is possible to consider instead of single
$\del$  and $tp_{k,\del}$  finite sequences of
ordinals $\del_1\nek\del_j$.  The corresponding
definitions of $\calL_{k,\del_1\nek \del_j}$,
$\gra_{k,\del_1\cdots \del_j}$  and
$tp_{k,\del_1\cdots\del_j}(\xi)$  are obvious.
The next lemma is similar to Lemma 2.1.

\proclaim Lemma 2.2.  Suppose that $j,k\le n$
and $k-j>2$. Let $\langle \del_1\nek\del_j\rangle$,
$\langle \gam_1\nek \gam_j\rangle$  be so that 
$$tp_k(\gam_1)=tp_k(\del_1)\leqno(a)$$
$$\hbox{for every}\ i<j\leqno(b)$$
$$tp_{k-i,\langle\gam_1\nek\gam_i\rangle}
(\gam_{i+1})= tp_{k-i,\langle\del_1\nek\del_i\rangle}
(\del_{i+1})\ .$$
Then for every $\xi <\lam$  there is $\rho<\lam$
such that
$$tp_{{k-j,}_{\langle\gam_1\nek\gam_j\rangle}}
(\rho)=tp_{{k-j,}_{\langle\del_1\nek\del_j\rangle}}
(\xi)\ .$$

\proclaim Lemma 2.3. There are $\xi,\rho \in
(\kap_n^{+n+1},\kap_n^{n+2})$ such
that $tp_n(\xi)=tp_n(\rho)$.

\pr Easy, since $2^{\kap_n}=\kap_{n+1}$ and each
type can be viewed as a subset of $\kap_n$.\hfill$\square$

Now let us deal with all $n$'s simultaneously.

\subheading{Definition 2.4}  Let $\oxi=\langle\xi_n\mid n
<\ome\rangle$.  $\orho =\langle\rho_n\mid n<\ome\rangle$
be two sequences in $\prodl_{n<\ome}\kap_n^{n+2}$
and $n^*<\ome$.  We say that $\oxi$  and
$\orho$  are equivalent starting with $n^*$  iff
for every $i<\ome$  $tp_{3+i}(\xi_{n^*+i})=tp_{3+i}
(\rho_{n^*+i})$.  $\oxi$  and $\orho$  are
called equivalent if there is $n^*<\ome$ such that
$\oxi$ and $\orho$  are equivalent starting with
$n^*$.  

\subheading{Definition 2.5}  Suppose that $\oxi$,
$\orho$ are equivalent sequences and
$\ogam,\odel\in\prodl_{n<\ome}\kap_n^{n+2}$.
Then $\ogam,\odel$ are called equivalent over
$\oxi$, $\orho$  starting with $n^*<\ome$  iff
for every $i<\ome$
$$tp_{3+i,\xi_{n^*+i}}(\gam_{n^*+i})=tp_{3+i,
\rho_{n^*+i}}(\del_{n^*+i})\ .$$

\proclaim Lemma 2.6.  Suppose that $\oxi$,
$\orho$,  $\ogam$ $\in\prod_{n<\ome}\kap_n^{n+2}$
and $\oxi,\orho$  are equivalent.  Then there is
$\odel\in\prod_{n<\ome} \kap_n^{n+2}$  such that
$\ogam,\odel$  are equivalent over $\oxi$,
$\orho$. 

\pr Follows from Lemma 2.1.\hfill$\square$

\sect{3.~~The preparation forcing}

Suppose $\kap_0<\kap_1<\cdots<\kap_n<\cdots$,  
$\kap_\ome=\bigcup_{n<\ome}\kap_n\quad ,\quad
o(\kap_n)=\kap_n^{+n+3}$.
We define a forcing of the type used in [Git-Mag1],
[Git-Mag2].  This remains especially those of [Git-Mag2],
section 3, but it will be different in a few
important points which will be used in further 
construction.

For every $n<\ome$, let us fix a nice system
$\UU_n=\ll \calU_{n,\alp}\mid\alp <\kap_n^{+n+2}\rangle$,
$< \pi_{n,\alp,\bet}\mid \alp,\bet <\kap_n^{+n+2},
\calU_{n,\alp}\triangleleft \calU_{n,\bet}\gg$.
We refer to [Git-Mag1] for the definitions.

Let $n<\ome$  be fixed.  We define first a
forcing $Q_n$  for adding one element Prikry
sequences for $\calU_{n,\alp}$'s $(\alp
<\kap_n^{+n+2})$.

\subheading{Definition 3.1} A set of forcing
conditions $Q_n$  consists of all elements $p$
of the form
$$\langle\{ <\gam, p^\gam >\mid \gam
<\del\},g,T\rangle$$
where 
\item{(1)} $g\subseteq\kap_n^{+n+2}$ of
cardinality $<\kap_n$  
\item{(2)} $\del <\kap_n^{+n+2}$
\item{(3)} if $g\not=\emptyset$, then $g$  has a
$\le_n$-maximal element and $0\in g$.
Further we shall denote $g$  by $\supp (p)$,
the maximal element of $g$  by $mc(p)$  and $\del$
by $\del(p)$.
\item{(4)} $\del \ge \cup g$
\item{(5)} for every $\gam\in g$ $\ p^\gam
=\emptyset$
\item{(6)} for every $\gam\in\del \bks g$
$p^\gam$ is an ordinal less than
$\kap_n^{+n+2}$
\item{(7)} $T\in U_{n,mc(p)}$

Denote $T$  by $T(p)$.

\subheading{Definition 3.2}  Let $p,q\in Q_n$.
We say that $p$  extends $q$ $(p\ge q)$  if  
\item{(1)} $\del (p)\ge \del (q)$
\item{(2)} for every $\gam\in\del
(q)\bks\supp(q)p^\gam=q^\gam$
\item{(3)} if $\supp (p)\not= \emptyset$  then
$\supp (p)\supseteq\supp (q)$ and
$\pi_{n,mc(p),mc(q)}{}\;\;\tagg T(p)\subseteq T(q)$.
\item{(4)} if $\supp (p)=\emptyset$  then
$p^{mc(q)}\in T(q)$  and for every $\nu\in\supp
(q)\ p^\nu=\pi_{n,mc(q),\nu}(p^{mc(q)})$.

\subheading{Definition 3.3}  Let $p,q\in Q_n$.
We say that $p$  is a direct extension of
$q(p\ge^*q)$ iff 
\item{(1)} $p\ge q$
\item{(2)} $\supp (q)\not=\emptyset$  implies
$\supp (p)\not=\emptyset$.
\item{}

Intuitively, we are going to add a one element
Prikry sequence for every $\nu\in \supp (p)$.
But there will be $\kap_n^{+n+2}$ coordinates
where nothing happens.  Namely we are allowing
$p^\gam$'s to be ordinals even above $\kap_n$
and this intended to prevent putting Prikry
sequence over $\gam$.

The proofs of [Git-Mag1] generalize directly to the
present situation.  Notice that in the argument
with elementary submodel $N$ there, $p\in Q_n$
will not be contained in $N$  anymore. But the
argument still works since $p\in N$  implies
$\del (p)\in N$  and so $\sup (N)$  can be added
to the support of $p$.

\proclaim Lemma 3.4.  $\le,\le^*$ are partial
orders on $Q_n$.

\proclaim Lemma 3.5.  For every $\gam
<\kap_n^{+n+2}$ $p\in Q_n$ with $\supp (p)\not=\emptyset$
there are $\gam'>\gam$  and $q\ge^* p$  with
$\gam'\in\supp (q)$.

\proclaim Lemma 3.6.  $Q_n$  has a dense set
isomorphic to the forcing of Cohen subset to
$\kap_n^{+n+2}$.

\proclaim Lemma 3.7.  $\langle Q_n,\le^*\rangle$
is $\kap_n$-closed.

\proclaim Lemma 3.8.  $\langle Q_n,\le,\le^*\rangle$  
satisfies the Prikry condition.

Let us now put all $Q_n$'s  together. 

\subheading{Definition 3.9} A set of forcing
conditions $\calP$  consists of all elements $p$
of the form
$$\langle p_n\mid n<\ome\rangle$$
so that
\item{(1)} $p_n\in Q_n$ $(n<\ome)$
\item{(2)} there is $\ell<\ome$  such that for
every $n\ge\ell$  $\supp (p_n)\not=\emptyset$. 

Further let us call the least such  $\ell$ the length
of $p$ and denote it by $\ell (p)$.  

The definition of the ordering $\le,\le^*$ on
$\calP$  is the usual definition of the order of
a product.

The proofs of the following lemmas are direct
generalizations of those in [Git-Mag1]. 

\proclaim Lemma 3.10. For every $f\in\prodl_{n<\ome}
\kap_n^{+n+2}$ and $p\in\calP$  there are
$g\in\prodl_{n<\ome}\kap_n^{+n+2}$  $g>f$  and
$q\in\calP$  s.t. $q{}^{*}\!\!\ge p$ and for every
$n\ge \ell$  $g(n)\in\supp (p_n)$.

\proclaim Lemma 3.11.  $\langle\calP,\le\rangle$
satisfies $\kap_\ome^{++}$-c.c.

\proclaim Lemma 3.12.
$\langle\calP,\le^*\rangle$ is $\kap_0$-closed
and for every $n<\ome, p\in\calP$  with $\ell
(p)\ge n$ $\langle \{ q\in\calP\mid q\ge
p\},\le^*\rangle$  is $\kap_n$-closed.

\proclaim Lemma 3.13.  $\langle\calP,\le,\le^*\rangle$ 
satisfies the Prikry condition.

Let $G\subset\calP$ be generic.  For
$f\in\prod_{n<\ome} \kap_n^{+n+2}$,  $f\in V$
let us denote by $G^f=\langle p^{f(n)}\mid n<\ome$,
$p\in G\rangle$.

\proclaim Lemma 3.14.  For every $g\in (\prod_{n<\ome}
\kap_n^{+n+2})\cap V$  there is $f\in (\prod_{n<\ome}
\kap_n^{+n+2})\cap V$  $f>g$  such that $G^f$  is a
Prikry sequence.

The proof follows from Lemma 3.10.

\proclaim Lemma 3.15.  If $f_1,f_2\in
(\prod_{n<\ome}\kap_n^{+n+2})\cap V$  $f_1<f_2$  and
$G^{f_1},G^{f_2}$  are both Prikry sequences,
then for all but finitely many $n$'s
$G^{f_1}(n)<G^{f_2}(n)$.

\pr There is a condition $p\in G$  such that for
every $n\ge \ell (p)$ $f_1(n),f_2(n)\in\supp
(p)$.  The rest follows from the definition of
nice sequence.\hfill$\square$

\proclaim Lemma 3.16.  For every $f\in
(\prod_{n<\ome}\kap_n^{+n+2})\cap V$  there is a Prikry
sequence.

\pr Pick $f^*>f$  such that $G^{f^*}$  is Prikry
and for all but finitely many $n$'s $f^*(n){}_n\!\!\ge
f(n)$.  Then $\langle
\pi_{n,f^*(n),f(n)}(G^{f^*}(n))\mid
n<\ome\rangle$ will be such a sequence.\hfill$\square$ 

\sect{4.~~Projection of $Q_n$}

Our aim is to define a projection of $\calP$  to
a less rigid p.o. set which still produces all the
Prikry sequences that $\calP$  does.

Let us split $\calP$  back and deal first with
$Q_n$'s $(n<\ome)$.

Fix $n<\ome$.  For a while we'll drop the lower
index $n$, i.e. $Q_n=Q$, $\kap_n=\kap$  etc.

Pick some $\alp_0\not=\alp_1<\kap^{n+2}$  realizing
the same type over $\kap^{+n}$  in sense of Section
1.  The automorphisms we are going to construct
will move the condition
$$\langle\{<0,<\gg,<\alp_0,<\gg\}\ ,\quad\{0,\alp_0\},
\kap\rangle$$
to the condition
$$\langle\{ <0,<\gg,<\alp_1<\gg\},\quad\{0,\alp_1\},
\kap\rangle\ .$$ 
Since $U_{\alp_0}=U_{\alp_1}$,  we can extend
this to conditions of the form
$$\langle \{ <0,<\gg, <\alp_i,<\gg\}\ ,\ \{
0,\alp_i\}\ ,\ A\rangle$$
for every $A\in U_{\alp_i}$ $(i<2)$.  Also let
us move this way one element Prikry sequences,
i.e.  $\langle \{<0,p^0>, <\alp_0,p^{\alp_0}>\}\rangle$
maps on $\langle \{< 0,p^0>,
<\alp_1,p^{\alp_0}>\}\rangle$.

The basic idea for further extensions is as
follows.  Let $\bet_0<\kap^{+n+2}$  be above
$\alp_0$  in the projection ordering $\le_n$.
Consider the extension of $\langle \{ <0,<\gg,
<\alp_0,<\gg\}, \{ 0,\alp_0\},\kap\rangle$
it generates i.e.
$$\langle \{ <0,<\gg, <\alp_0,<\gg, <\bet_0,
<\gg\}, \{ 0,\alp_0, \bet_0\},\kap\rangle\ .$$

We would like to pick some $\bet_1$  similar to
$\bet_0$  and to map the last condition onto
$$\langle \{ <0,<\gg, <\alp_1,<\gg, <\bet_1
<\gg\}, \{ 0,\alp_1, \bet_1\},\kap\rangle\ .$$

There are two problems with this.  Firstly there 
may be no $\bet_1$  realizing the same type over
$\alp_1$  as $\bet_0$  does over $\alp_0$.
The second is that there may be lots of
$\bet_1$'s which are good candidates.  In this
case we should be careful not to pick two of
them $\bet'_1, \bet\tagg_1$  such that they have
a common part which is mapped in different
fashions to $\bet_0$.  This eventually may hide
all Prikry sequences and leave us with a trivial
forcing. 

We shall inductively define an equivalence
relation $\lft$ on a dense subset of $Q$  so that for
every $p,q\in Q$ if $p\longleftrightarrow q$
then for any $p'\ge p$  there are $p\tagg\ge p'$
and $q\tagg\ge q$ such that $p\tagg\lft
q\tagg$. Then we define a preorder $\to$ over
$Q$ which will combine $\le$  and
$\longleftrightarrow$.  In particular if   $p\le q$ or
there is $q'\le q$ $p\longleftrightarrow q'$, then 
$p\to q$.

Let $\{\langle p^0_i,p^1_i,p^2_i\rangle \mid
i<\kap^{+n}\}$ be an enumeration with unbounded
repetitions of triples of $Q$ so that
$$p_0^j=\langle \{<0,<\gg, <\alp_j,<\gg\},\{0,\alp_j\},
\kap\rangle$$
for $j<2$  and $p^2_0=p^0_0$.  

Define $\longleftrightarrow$  by induction on
$i<\kap^{+n}$.  For $i=0$  let $p^0_0\lft p^1_0$. 
Suppose that for every $i'<i$  the
triple $\langle p^0_{i'}p^1_{i'}p^2_{i'}\rangle$
was considered.  If either (1) not
$p^0_i\longleftrightarrow p^1_i$ or (2)
$p^0_i\longleftrightarrow p^1_i$  and
$p^2_i\not\ge p^0_i$  then nothing new happened
and we continue to the next triple.  Suppose
otherwise.  Then $p^0_i\longleftrightarrow
p^1_i$,  $p^2_i\ge p^0_i$.

\subheading{Case 1} $\supp (p^0_i)=\supp
(p^2_i)$.

Let $\del =\del (p_i^2)+\del (p^1_i)+1$.
Pick $p,q\in Q$  such that 
\item{(a)} $p\ge p^2_i$
\item{(b)} $\supp (p)=\supp (p^2_i)$
\item{(c)} $q\ge p^1_i$
\item{(d)} $\supp (q)=\supp (p^1_i)$
\item{(e)} $\del (p)=\del (q)\ge \del +i$
\item{(f)} for every $\gam\in\del (p)\bks \del
(p^2_i)$, $\xi\in\del (q)\bks\del (p^1_i)$.

$q^\xi,p^\gam  >\max (\kap, \{ t^\gam\mid t\in
Q\quad \hbox{appeared in the same triple before stage}
\quad i\quad {\rm and}\quad\gam <\del (t)\}\quad
{\rm and}\quad q^\xi\not= p^\gam$.

Set $p\longleftrightarrow q$.

\subheading{Case 2} $\supp(p^0_i)\not=\supp
(p^2_i)$. 

\subheading{Subcase 2.1}  $\supp
(p^2_k)=\emptyset$ and $\del (p^0_i)=\del (p^1_i)$.

Then $p^2_i$  is obtained from $p^0_i$  by
adding a one element sequence to $mc(p^0_i)$
and projecting it to the rest of supp $(p_i^0)$.
Add this sequence to $mc(p^1_i)$ and project it
to the rest of $\supp (p^1_i)$.   
Let $q$  be the condition obtained from $p^1_i$
this way.  Set $p^2_i\longleftrightarrow q$.
\subheading{Subcase 2.2} $\supp (p^2_i)=\emptyset$
and $\del (p^0_i)<\del (p^2_i)$. 

Then we add $(p^2_i)^{mc(p^0_i)}$ to $p^1_i$ as
a one element sequence at coordinate $mc(p^1_i)$  and
project it to all the rest of $\supp(p^1_i)$.
Let $q'$ be a condition obtained from $p^1_i$
using this process.

Now pick $p\ge p^1_i$  and $q\ge q'$  as in Case
1.  Set $p\longleftrightarrow q$.

\subheading{Subcase 2.3}  $\supp
(p^2_i)\supsetneqq \supp (p^0_i)$.

Consider the type of $mc(p^2_i)$ over
$mc(p^0_i)$. Let $\del =\kap +i+\sup \{ t^\gam +\del
(t)|t\in Q,\gam <\del (t)\hbox{and}\ t\ \hbox{or some
condition which was defined to be equivalent
to}\ t\ \hbox{appeared before or at}$\hb
$\hbox{the stage}\ i\}$. 

If there is some $\tau >\del$  realizing the
same type over $mc(p^1_i)$ as $mc(p^2_i)$  over
$mc(p^0_i)$,  then proceed as follows.  Add $\tau$
to the support of $p^1_i$  and for every
$\gam,\del (p^1_i)\le\gam <\tau$  add
$\langle\del\rangle$  to coordinate $\gam$.  Let
$q$  be a condition obtained this way.  Let
$p=p^2_i\cup \{ \langle \gam, \del\rangle \mid
\del (p^2_i)\le \gam\le\tau\}$. Set
$p\longleftrightarrow q$.

Assume now that no $\tau >\del$  realizes the same
type over $mc(p^1_i)$  as $mc(p^2_i)$  does over
$mc(p^0_i)$. Let $p=p^2_i\cup \{\langle \gam ,\del\rangle
\mid \del(p^2_i)\le\gam <\del\}$ and
$q=p^1_i\cup \{ \langle \gam, \del \rangle
\mid\del (p^1_i)\le\gam <\del \}$.  Consider $T(p)$
i.e. the set of $\calU_{\kap ,mc(p)}$  measure
one of $p$.  For every $\nu\in T(p)$, denote by
$p^\cap \langle\nu\rangle$  the extension $p$  obtained
by adding $\nu$ to the maximal coordinate of $p$  and
projecting it to every coordinate of $\supp
(p)$.  We also consider $q^\cap(\pi_{mc(p),mc(p^0_i)}
(\nu))$ which is defined similarly.  Set
$$p^\cap\langle\nu\rangle\cup\{ \langle\del,\del \rangle\}
\longleftrightarrow q^\cap(\pi_{mc(p),mc(p^0_i)}(\nu))\cup
\{ \langle \del,\nu\rangle\}\ .$$
This completes the definition of $\longleftrightarrow$.

Notice that always $p\longleftrightarrow q$  implies
$\del (p)=\del (q)$ 

Let $D=\{ p\in Q\mid\exists q\in Q\ p\longleftrightarrow
q\}$. 

\proclaim Lemma 4.1.  $D$  is dense in $\langle
Q,\le\rangle$.

\pr Let $p\in Q$.  Extend it to some $t$
stronger than $p^0_0$.  Then the triple $\langle
p^0_0,p^1_0,t\rangle$  appears at some stage $i$
of the construction.  Then either for some $t'$
$t\longleftrightarrow t'$ or for some $t^*\ge t$
$\ t^*\longleftrightarrow t'$.\hfill$\square$

\proclaim Lemma 4.2.  For every $p\in D$. If $q\not=
q'\in D$  are such that $p\lft q$, $p\lft q'$.
Then $q,q'$ are incompatible. 

\pr Suppose otherwise.  Then let
$i<j<\kap^{+n}\ q_i,q_j$  be such that $p\lft
q_i=q$,  $p\lft q_j=q'$, $q_i$  appeared at stage
$i$  and $q_j$  at stage $j$.  Since $i<j$  $p$
is supposed to be really stronger or
incompatible with every element appeared before
stage $j$.  But this includes $p$  itself.
Which is impossible.  Contradiction.\hfill$\square$   

\proclaim Lemma 4.3.  Let $p,q\in D$  $p\lft q$
and $p\not= p^j_0$  for $j< 2$. Then $p,q$  are
incompatible. 

The proof follows from the definition of $\lft$.

Let $D^*$  denote $D/\lft$.  Define a preorder
$\to$  over $D^*$  and actually over $Q$.

\subheading{Definition 4.4} $p\to q$  iff there
exists $n<\ome$ and a sequence $\langle p_k\mid
k<n\rangle$ so that 
\item{(1)} $p_0=p$
\item{(2)} $p_{n-1}=q$
\item{(3)} for every $k<n-1$
$$p_k\le p_{k+1}\quad \hbox{or}\quad p_k\lft
p_{k+1}\ .$$ 
See diagram:
$$\matrix{p_{n-2}&\lft&p_{n-1}\cr
\vee|&&\cr
p_{n-3}&\lft&p_{n-4}\cr
\cdots&\cdots&\cdots\cr
\cdots&\cdots&\vee|\cr
p_4&\lft&p_5\cr
\vee|&&\cr
p_3&\lft&p_2\cr
&&\vee|\cr
p_0&\lft&p_1\cr}$$
Clearly $\to$  is reflexive and transitive.

\proclaim Lemma 4.5. $\langle D^*,\to\rangle$
is a nice suborder of $\langle Q,\le \rangle$,
i.e. for every generic $G\subseteq Q$
$\ \{p/\lft\mid p\in G\cap D\}$ is a generic for
$\langle D^*,\to \rangle$.

\pr It is enough to show the following: 

For every $E^*\subseteq D^*$  dense open in
$\langle D^*,\to \rangle$ $E=\{ p\in Q\mid
p/\lft \in E^*\}$  is dense in $\langle Q,\le
\rangle$.  So let $E^*\subseteq D^*$  be a dense
open set and $s\in Q$.  First we extend $s$  to
some $p\in D$. Then there is $q\in D$ such that
$p\to q$  and $q/\longleftrightarrow \in E^*$. Let
$\langle p_k\mid k<n\rangle$  be as in Definition 4.4. 
Assume for simplicity that $n=\infty$.  At
the first step find $p'_6\ge p_6$  and $p'_4\ge
p_4$  so that $p'_6\longleftrightarrow p'_4$.
This is possible since the triple $\langle
p_5,p_4,p_6\rangle$  appears at some stage of
the definition of $\longleftrightarrow$.  See
the diagram:
$$\matrix{&&&&&&p'_6&&&&&\cr
&&&&&&\vee|\cr
&&q&=&p_7&\longleftrightarrow&p_6&&&&&\cr
&&&&&&\vee|\cr
p\tagg_4&\ge&p'_4&\ge&p_4&\longleftrightarrow
&p_5&&&&\cr
&&&&\vee|&&&&&&&\cr
&&&&p_3&\longleftrightarrow&p_2&\le&p_2\tagg&\le&p_2
\taggg\cr 
&&&&&&\vee|&&&&&\cr
p\taggg_0&\ge&p&=&p_0&\longleftrightarrow&p_1&&&&&\cr}$$

Then we shall find $p\tagg_4\ge p'_4$  and $p\tagg_2\ge
p_2$  such that $p\tagg_2\longleftrightarrow
p\tagg_4$.  They exist by the same reason.  And
finally we pick equivalent $p\taggg_0\ge p_0=p$  and
$p\taggg_2\ge p\tagg_2$. Now, as $q/\lft \in E^*$, 
$p_6\lft q$ and $p'_6\ge p_6$,  $p'_6/\lft\in E^*$.  Then
$p'_4/\lft \in E^*$  and therefore $p\tagg_4/\lft$. 
This implies $p\tagg_2/\lft\in E^*$  and then
$p\taggg_2/\lft\in E^*$.  But then $p\taggg_0/\lft
\in E^*$.  Therefore $p\taggg_0\in E$  and we
are done.\hfill$\square$

The same proof gives the following.

\proclaim Lemma 4.6.  For every $p\in Q$ the set
$\{ q\mid p\le q$  or $p,q$  are $\to$ incompatible
$\}$ and is dense in $\langle Q,\le\rangle$. 

Let $p\in Q$  be a condition with $\supp (p)\not=
\emptyset$. For $\nu\in T(p)$, (recall that
$T(p)$  is a set of measure one corresponding to
the maximal coordinate of $p$) let us denote by
$p^\cap \langle\nu\rangle$  the condition
obtained from $p$  by adding $\nu$  to the
maximal coordinate of $p$  and projecting it to
all permitted for $\nu$  coordinates in $\supp (p)$.

Suppose $p,q\in Q$,  $\supp p\not=\emptyset$,
$\supp q \not=\emptyset$  and $p\longrightarrow
q$.  Let $\vec p=\langle p_k\mid k<n\rangle$  be
a sequence witnessing $p\longrightarrow q$.  Let
$\gam\in\supp p$.  We define the orbit of $\gam$
via $\vec p$  as follows:  $\gam_0=\gam$, $\gam_1$
the coordinate of $p_1$  corresponding to
$\gam_0$, for every $0<k<n$  let $\gam_k$  be
the coordinate corresponding to $\gam_{k-1}$  in
$p_k$.  Let us call $\gam_{n-1}$ the $\vec p$-image
of $\gam$  in $q$.

The principal point will be to show that $\vec p$-image
of $\gam$  in $q$  actually does not depend on a
sequence $\vec p$ witnessing $p\longrightarrow
q$.  This implies, in particular that if $p\le 
q$, then the only image of $\gam$  in $q$  is
$\gam$  itself.  Which in turn implies that Prikry
sequence are produced by $\langle Q,\longrightarrow
\rangle$, since then for different $\nu_1,\nu_2$
$p^\cap \langle \nu_1\rangle$  and
$p^\cap\langle \nu_2\rangle$  are incompatible
in $\langle Q,\longrightarrow\rangle$.

\proclaim Lemma 4.7.  Let $p\longrightarrow q$
$\supp p, \supp q\not=\emptyset$  and $\gam\in\supp
p$.  Then the image of $\gam$  in $q$ does not
depend on a particular sequence witnessing
$p\longrightarrow q$.

\pr Suppose otherwise.  Let us pick $q$  with
$\del (q)$  as small as possible so that
$p\longrightarrow q$,  $\supp (q)\not=\emptyset$
and for some $\gam\in\supp (p)$  image of $\gam$  in
$q$ depends on particular sequence witnessing
$p\longrightarrow q$.  Pick two sequences of the
least possible lengths  $\vec s=\langle
s_i\mid i<s^*<\ome\rangle$  and $\vec t=\langle
t_i\mid i<t^*<\ome\rangle$  witnessing
$p\longrightarrow q$  and such that $\vec
s$-image of $\gam$ in $q$  is different from $\vec
t$-image of $\gam$  in $q$.     

W.l. of $g$.  $q=s_{s^*-1}>s_{s^*-2}$  or
$q=t_{t^*-1}>t_{t^*-2}$.  Since, otherwise,
$q\longleftrightarrow s_{s^*-2}$  and
$q\longleftrightarrow t_{t^*-2}$.  But by the
construction of $\longleftrightarrow$ the above
implies $s_{s^*-2}=t_{t^*-2}$  and we can then
replace $q$  by $s_{s^*-2}$  reducing the
lengths of the witnessing sequences.  Notices
that it is impossible to have both $q>s_{s^*-2}$
and $q>t_{t^*-2}$.  Since then by the definition
of $Q$ $s_{s^*-2}\ge t_{t^*-2}$  or $t_{t^*-2}$
$\ge s_{s^*-2}$  (we are ignoring sets of
measure).  Let $s_{s^*-2}\ge t_{t^*-2}$.  Then
$p\longleftrightarrow s_{s^*-2}\gam$  has different
images in $s_{s^*-2}$  and $\del
(s_{s^*-2})<\del (q)$.  So let us assume that
$q>s_{s^*-2}$  and $q\longleftrightarrow
t_{t^*-2}$.  Denote $s_{s^*-2}$  by $s$  and
$t_{t^*-2}$  by $t$.  Consider the stage in the
construction of $\longleftrightarrow$  where $q$
and $t$  were made equivalent.  Suppose it was
at a stage $i$.  So the triple $\langle
p_{i0},p_{i1},p_{i2}\rangle$  were considered.
Since $\del (q)>\del (s_k)$,  $\del (t_\ell)$
for every $k\le s^*-2,\ \ell\le \ell^*-3$,  the
equivalences inside the sequences $\langle s_k\mid
k<s^*-1\rangle$  and $\langle t_\ell \mid \ell
<t^*-2\rangle$  were obtained on earlier
stages. 

We claim that for some $j<2$  $p_{ij}\ge s$
(once again we are ignoring sets of measure one
in the condition considered).  Since otherwise
the definition of $\longleftrightarrow$  will
make $q$ incompatible with $s$.  By the same
reason for $j<2$  $p_{ij}\ge t_{t^*-3}$.  But
$\gam$  is already embedded into $s$  and
$t_{t^*-3}$  it is embedded into $p_{i0}$  and
$p_{i1}$.  Also $p_{00}\longleftrightarrow p_{01}$
and $\del (p_{00})=\del(p_{01})<\del (q)$.  Hence
picking $p_{00}$  or $p_{01}$  instead $q$  we
obtain a contradiction to the minimality of $\del
(q)$.\hfill $\square$

\sect{5.~~The projection of $\calP$}

Basically we are going to now put together the
projections defined in the previous section for
each $n<\ome$. 

Add the lower index $n$  to everything defined
in Section 3, i.e. $\alp_{0n},\alp_{1n},D_n,D_n^*$,
$\lft_n,\longrightarrow_n$.  Set $D=\{
p\in\calP|$  for every $n$ $p\rhookup Q_n\in
D_n\}$.  Denote $p\rhookup Q_n$  by $p_n$.

Let $p,q\in Q$  and $p\lft_nq$.  Suppose that
$\supp (p)\not=\emptyset$.  Then for some
$k<\ome$  $mc(p_n), mc(q_n)$  are realizing the
same type over $\kap_n^{+k+2}$.  We denote this
by $p\lft_{n,k} q$.

\subheading{Definition 5.1}  Let $p,q\in D$.
Then $p\lft q$  iff
\item{(1)} for every $n<\ome$ $p_n\lft_n q_n$
\item{(2)} there is a nondecreasing sequence
$\langle k_n\mid n<\ome \rangle$,
$\lim_{n\to\infty}k_n=\infty$  such that for
every $n\ge \ell (p)$ $p_n\lft_{n,k_n}q_n$.

Clearly, $\lft$ is an equivalence relation on
$D$.  Set $D^*=D/\lft$.  Let us now turn to the
ordering $\to$ on $D^*$.  For $p,q\in Q_n$
$p\to_nq$  with $\supp (q)\not=\emptyset$  let
us denote by $p\to_{n,k}q$  the fact of
existence of a sequence witnessing $p\to_nq$
such that the equivalent members are at least
$\lft_{n,k}$  equivalent.

\subheading{Definition 5.2}  
Let $p,q\in D$.  Then $p\to q$ iff there is a
finite sequence of conditions $\langle r_k\mid
k<m<\ome\rangle$  so that 
\item{(1)} $r_0=p$
\item{(2)} $r_{m-1}=q$
\item{(3)} for every $k<m-1$ $r_k\le r_{k+1}$
or $r_k\longleftrightarrow r_{k+1}$,
or schematically:
$$\matrix{r_{m-2}&\longleftrightarrow&r_{m-1}=q\cr
\vee |&&\cr
r_{m-3}&\longleftrightarrow&r_{m-4}\cr
&\cdots&\cr
&&\vee|\cr
r_4&\longleftrightarrow&r_5\cr
\vee |&&\cr
r_3&\longleftrightarrow&r_2\cr
&&\vee|\cr
p=r_0&\longleftrightarrow&r_1\cr}$$

\proclaim Lemma 5.2.1.  Let $p,q\in D$.
If $p\to q$  then
\item{(1)} for every $n<\ome$  $p_n\to_nq_n$ 
\item{(2)} there is a nondecreasing sequence
$\langle k_n\mid n<\ome \rangle$, $\lim_{n\to\infty}
k_n=\infty$  such that for every $n\ge \ell(q)$ 
$p_n\to_{n,k_n}q_n$.

\pr Let $\langle r_k\mid k<m\rangle$  be a
sequence witnessing $p\longrightarrow q$.  Then
for every $k<m-1$  $r_k\le r_{k+1}$  or
$r_k\longleftrightarrow r_{k+1}$. In the latter
case there exists a nondecreasing sequence
$\langle\ell_{kn}\mid n<\ome\rangle$
$\lim_{h\to\infty}\ell_{k,n}=\infty$  such that
for every $n\ge\ell (r_k)\
r_{kn}\longleftrightarrow_{n\ell_{kn}}r_{k+1n}$.
Set $\ell_n=\min\{\ell_{k,n}\mid k<m\}$ for
every $n<\ome$.  Then $\langle\ell_n\mid n<\ome\rangle$
is a nondecreasing sequence converging to
infinity.  Also, for every $n\ge \ell (q)$
$p_n\longrightarrow_{n\ell_n}q_n$.  Since
$\langle r_{kn}\mid k<m\rangle$  be a witnessing
sequence for this.\hfill$\square$   

\proclaim Lemma 5.2.2.  $D$  is a dense subset of
$\calP$.

\pr Let $p\in\calP$.  Triples $\langle
p^0_{0,n},p^1_{0,n},p^2_n\rangle$  appear in
enumerations of $Q_n$'s $(n<\ome)$.  So for each $n$
there will be $p'_n\ge_np_n$  in $D_n$.
So the natural candidate for $q\ge p$ in $D$  is
$\langle p'_n\mid n<\ome\rangle$.  The only
problem is that by the definition of $\calP$,
starting some $\ell <\ome$  the supports of
$p'_n$'s are supposed to be nonempty.  

Consider $\ell(p)$  (the least $\ell <\ome$
s.t. $\supp (p_n)\not=\emptyset$  for every $n\ge\ell$).
Let $n\ge \ell (p)+2$.  Since
$\alp_{0,n},\alp_{1,n}$  are realizing the same
type over $\kap_n^{+n}$, there will be
arbitrarily large (below $\kap_n^{+n+2}$)
$\bet$'s realizing over $\alp_{1,n}$  the same
$\kap_n^{+n-1}$-type as $mc(p_n)$  does over
$\alp_{0,n}$.  So, we are in the situation of
Subcase 2.2. of the definition $\lft_n$.  Hence
$p'_n\ge p_n$  with $\supp_np'_n=\supp p_n$  is
picked and added to $D_n$.  This means that $q$
consisting of such $p'_n$ will be as desired.
\hfill$\square$

\proclaim Lemma 5.3.  $\langle D^*,\to\rangle$
is a nice suborder of $\langle\calP\le\rangle$.

\pr We proceed as in Lemma 4.5.  Let
$E^*\subseteq D^*$  be a dense open subset of
$\langle D^*,\to\rangle$.  We need to show that
$$E=\{p\in\calP \mid p/\lft \in E^*\}$$
is dense in $\langle \calP,\le\rangle$.  Let
$s\in\calP$.  Extend it to some $p\in D$. Then
for some $q\in D$  $p\to q$  and $q/\lft \in
E^*$.  For $n$'s below $\ell(q)$ we just repeat
the argument of Lemma 4.5.  For $n$'s above
$\ell (q)$  pick a sequence $\langle k_n\mid
n<\ome\rangle$  as in 5.2.1(2). So $p_n\to_{n,k_n}q_n$
for every $n\ge \ell (q)$.  Consider the diagram
like those of Lemma 4.5.
$$\matrix{&&&&&&h'&&&&\cr
&&&&&&\vee|&&&&\cr
&&&&q_n&\buildrul n,k\under\lft&h&&&&\cr 
&&&&&&\vee|&&&&\cr
f\tagg&\ge&f'&\ge&f&\buildrul n,k\under\lft&g&&&&\cr
&&&&\vee|&&&&&&\cr
&&&&d&\buildrul n,k\under\lft&c&\le&c\tagg&\le&c\taggg\cr
&&&&&&\vee|&&&&\cr
a\taggg&\ge&p_n&=&a&\buildrul
n,k\under\lft&b&&&&\cr}$$

We can find $f'\buildrul n,k-1\under\lft h'$
since $g \buildrul n,k\under\lft f$ and $h\ge
g$.  So there are arbitrarily large
$\xi<\kap_n^{+n+2}$  realizing the same
$\kap_n^{+(k+2)-1}$ -- type over $mc(f)$  as
$mc(h)$  realizes over $mc(g)$.  Applying this
to the triple $\langle g,f,h\rangle$ we obtain
such $f',h'$.  Now the same argument for $\langle
d,c,f'\rangle$  produces
$f\tagg\buildrul n,k-1\under\lft c\tagg$ and, finally,
we reach
$a\taggg\buildrul n,k-1\under\lft c\taggg$ with
$a\taggg\ge p_n$.    

Notice that $q_n\buildrul n,k-1\under\longrightarrow a
\taggg$.  Since 
$$\matrix{a\taggg&\buildrul n,k-1\under\lft&c\taggg\cr
&&\vee|\cr
f\tagg&\buildrul n,k-1\under\lft&c\tagg\cr
\vee|&&\cr
f'&\buildrul n,k-1\under\lft&h'\cr
&&\vee|\cr
q_n&\buildrul n,k\under\lft&h\cr}$$
Denote $a\taggg$  by $p^*_n$.

Doing this analyses for every $n$  we obtain a
condition $\langle p^*_n\mid n<\ome\rangle =p^*$
such that $p^*\in D,\ p^*\ge p$  and $q\longrightarrow
p^*$.  Hence $p^*/\longleftrightarrow\in E^*$
and so $p^*\in E$.\hfill$\square$ 	 

\sect{6.~~ Extensions of Elementary Submodels in
the Forcing With $\calP$} 

Let $N$  be an elementary submodel of a large
enough piece of the universe of cardinality at
most $\kap_\ome$.  Let $G$ be a generic subset
of $\langle \calP, \le \rangle$ or of the
projection of it defined in Section 4.  We like
to have $N[G]$  close to $N$.  More specific, we
like $N[G]$  to have the same ordinals as $N$
and to be closed under $\ome$-sequences.
Unfortunately, it is impossible to achieve.
Thus, if $N[G]$  is $\ome$-closed then $|N|<\kap_\ome$
and then $G$ will add new ordinals to $N$. 

The compromise will be to to deal with $N$'s of
cardinality $\kap_\ome$ together with its
elementary submodel $N^*$ so that   
\item{(1)} $\cup (N^*\cap \kap^+_\ome)=\cup
(N\cap\kap_\ome^+)$
\item{(2)} for every $n<\ome$ $\ \cup (N^*\cap
\kap_n^{+n+2})=\cup (N^*[G]\cap \kap_n^{+n+2})$
\item{(3)} ${}^\ome\!N^*[G]\subseteq N^*[G]$.
\item{}

Let us turn to the definition.

\subheading{Definition 6.1} A model $N$  is
called a good model iff
\item{(1)} $N\cap \kap^+_\ome$  is an ordinal.
\item{(2)} there exists an increasing continuous
sequence of submodels $\langle N_i\mid i\le\del
\rangle$  with $N_\del =N$  such that
\item{(2a)} $\langle N_i\mid i<j\rangle\in N_{j+1}$
for every $j<\del$
\item{(2b)} $\bigcupl_{i<\del}N_i=N$
\item{(2c)} $\kap_0>cf \del=\del >\aleph_0,\del^\ome=
\del$ 
\item{(3)} there exists an increasing continuous
sequence of submodels $\langle N_i^*\mid
i\le\del\rangle$  of $N$  (the same $\del$  as
in (1)) such that
\item{(3a)} $\langle N^*_i\mid i <j\rangle\in
N^*_{j+1}$
\item{(3b)} $\del =|N^*_i|$  for every $i<\del$
\item{(3c)} ${}^\ome\!N^*_{i+1}\subseteq N^*_{i+1}$
for every $i<\del$
\item{(3d)} $\{N_i\}\cup\del\subseteq
N^*_{i+1}\subseteq N_{i+2}$  for every $i<\del$
\item{(3e)} $\cup (N^*_\del
\cap\kap^+_\ome)=\cup (N\cap \kap^+_\ome)$.
\item{}

Further we call $\ll N_i\mid i<\del>$,  $\langle
N^*_i\mid i\le\del \gg$ a good sequence for $N$
and denote $N^*_\del$  simply by $N^*$.  Let us
refer, also, to $\langle N,N^*\rangle$ as to a
good pair.  For the main results in the next
section only $N^*$  will be important.	

For a generic $G\subseteq\calP$  and a model
$N$, we mean by $N[G]$  the set $\{\tau [G]\mid\tau\in
N$, $\tau$ is a name $\}$. 

\proclaim Lemma 6.2.  Le $D$  be a $*$-dense
open subset of $\calP$, (i.e. dense open in
$\langle\calP,\le^*\rangle$)  then for every
$n<\ome$  there are $D_n\subseteq
\calP\rhookup n+1$,  $D^n\subseteq\calP\bks n$
such that
\item{(a)} $D^n$  is $\kap_n$-weakly closed.
\item{(b)} $D\times D^n\subseteq D$
\item{(c)} $D_n\times D^n$  is $*$-dense in
$\calP$.
\item{}

The proof follows from Lemma 3.12, since
$|\calP\rhookup n+1|<\kap_{n+1}$.

\subheading{Definition 6.3} A condition $p\in\calP$
is called $*$-generic for a pair of models
$\langle N,N^*\rangle$  iff for every $\calP$-names
$\tau,\tau^*$  of ordinals such that $\tau\in
N$,  $\tau^*\in N^*$ $p$  forces in the forcing
$\langle\calP,\le^*\rangle$  ``$\tau [\buildrul^\sim\under
G]\in N$, $\tau^*[\buildrul^\sim\under G]\in
N^*$".

\proclaim Lemma 6.4.  Let $\langle N, N^*\rangle$
be a good pair.  Then there is a $*$-generic
condition for $\langle N,N^*\rangle$.

\pr Let $\langle N_i\mid i\le\del\rangle$,
$\langle N_i\mid i\le\del\rangle$  be sequences
witnessing the goodness of the pair $\langle
N,N^*\rangle$.  Let $\langle D_{i,j}\mid j<\kap_\ome
\rangle$, $\langle E_{i,j}\mid j<|N^*_i|\rangle$  be
enumerations of $*$-dense open sets of $N_i$ and
$N^*_i$  respectively.  For every$j$,
$\kap_n\le j<\kap_{n+1}$.  Using Lemma 6.3,
replace each $D_{i,j}$  by product $D'_{i,j}\times
D\tagg_{i,j}$  with $D\tagg_{i,j}$
$\kap_{n+2}$-weakly closed.

W.l. of $g$. let us assume that for every
$i<\del$  $\langle D_{i,j}\mid j<\kap_\ome\rangle$,
$\langle E_{i,j}\mid j<|N^*_i|\rangle$  belong
to $N^*_{i+1}$.

We define now by induction  a $*$-increasing
sequence of conditions $\langle p_i\mid i\le
\del\rangle$.  Let us describe only the first
stage of the construction, since each successor
stage will be like this and at limit stages the
union will be taken. 

Pick $q\in N^*_1$  such that for every $j<\kap_\ome,
\kap_n\le j<\kap_{n+1}$  $q$ belongs to
$D\tagg_{0j}$.  Let $p_1{}\;\ge^*q$, $p_1\in N_2$
be such that $p_1\in\bigcap\limits _{j<|N^*_0|}E_{0j}$.

Let us show that $p=p_\del$  is a $*$-generic
for $\langle N, N^*\rangle$.  There is no
problem with $N^*$, since it meets each
$*$-dense set of $N^*$.  Let $\tau$  be a name
of ordinal in $N$.  Then for some $i<\del$
$\tau\in N_i$.  Let $D=\{q\in \calP\mid q$ 
decides the value of $\tau$  in the forcing
$\langle\calP,\le^*\rangle\}$.  Then for some
$j<\kap_\ome$  $D\supseteq D'_{ij}\times D\tagg_{ij}$.
Pick $n<\ome$  so that $\kap_n\le j <\kap_{n+1}$.
By the choice of $p$, $p\ge p_{i+1}$  and $p_{i+1}$
is in $D\tagg_{ij}$.  Since $D'_{ij}\subseteq\calP\rhookup
n+1$,  $|D'_{ij}|<\kap_{n+1}$.  So
$D'_{ij}\subseteq N_i$.  Then $p$  forces in
$\langle \calP, \le^*\rangle$  $\tau
[\buildrul^\sim\under G]\in N_i$.\hfill$\square$ 

Now let us turn to the real forcing, i.e.
$\langle\calP,\le \rangle$.

\proclaim Lemma 6.5.  Let $N$  be a good model
and $p\in\calP$  be $*$-generic over $N$. Then
$p$  forces in $\langle \calP,\le \rangle$   
$$\tagg N[\buildrul^\sim\under G]\cap On=N\cap
On\tagg\ .$$

\pr Let $\tau\in N$  be a name of an ordinal.
By Lemma 3.13, for every $q\in\calP$  there is
$q'\ge^* q$  and $\ome >n\ge \ell (q')$  such
that the extension of $q'$ obtained by adding an
element of Prikry sequence of $mc(q'_n)$
already decides the value of $\tau$.  Certainly,
this value may depend on particular element added, 
but Prikry sequence consists of elements below
$\kap_\ome$  and in our case even below
$\kap_n$.  So, if $q'\in N$  then every such
extension is in $N$  as well.  Hence $\tau
[G]\in N$.

Let us show that it is possible to find such
$q'\in N$ below $p$.  Let $D$  be a dense open
set of conditions of this type.  Then, for some
$i<\del$  $D\in N_i$ where $\del$,  $N_i$  are
as in the definition of a good model.  Now split
$D$  into $D'_{ij}\times D\tagg_{ij}$  and
proceed as in Lemma 6.4.\hfill$\square$

Let us turn to $N^*$ now.  Suppose that $\langle N,
N^*\rangle$  is a good pair and $p\in\calP$  is
$*$-generic for $\langle N,N^*\rangle$.
Consider $\langle mc_n(p)\mid n\ge \ell
(p)\rangle$, the sequence of maximal coordinates
of $p$  and $\langle A_n\mid n\ge \ell
(p)\rangle$  the sequence of corresponding sets of 
measures one.  By the choice of $p$  and Lemma
3.13, for every $\calP$-name $\tau\in N^*$  of
an ordinal there is $n$,  $\ome >n\ge \ell (p)$
such that for every $\nu\in A_n$,  $p^\cap
\langle \nu\rangle$  decides $\tau$,  where
$p^\cap \langle \nu\rangle$  is the extension of
$p$  obtained by adding $\nu$  to $mc_n(p)$ and
projecting $\nu$ to coordinates in
$\supp_n(p)$.

For a generic $G\subseteq\calP$,  set $N^*[G]=\{\tau
[G]\mid\tau\in N^*\}$. 
Let $\langle N^*_i\mid i<\del\rangle$  be a
sequence witnessing goodness of $N^*$.  By Lemma
6.4, $p$  can be chosen to be a supremum of
$\langle p_i\mid i<\del\rangle$  so that each
$p_i$  is $*$-generic over $N_i^*$  and $p_i\in
N^*_{i+1}$.  Suppose that this is the case.  

\proclaim Lemma 6.6.  Let $\lam\in N^*$  be a
regular cardinal and suppose that $\lam\notin
\{\kap_n\mid n<\ome\}$.  Then for every generic
$G\subseteq \calP$ with $p\in G$ 
$$\cup (N^*[G]\cap\lam)=\cup (N^*\cap\lam)\ .$$

\pr Let $\tau\in N^*$  be a name of an ordinal
less than $\lam$.  Find $i<\lam$  such that
$\tau\in N^*_i$.  Since $p_i$  is $*$-generic
for $N^*_i$,  there will be $n$,  $\ell (p_i)\le
n<\ome$  and a set of measure one $A_{in}$  in
$p_i$  corresponding to $mc_n(p_i)$  such that
for every $\nu\in A_{in}$  $p_i^\cap
\langle\nu\rangle$ decides the value of $\tau$.
But $p_i\in N^*_{i+1}\subseteq N^*$.  So, there
will be some ordinal $\tau^*\in N^*\cap\lam$
which will bound all the possibilities for the
value of $\tau$.\hfill$\square$

\subheading{Remark 6.7}  The same argument works
for every limit $i<\del$.

\proclaim Lemma 6.8.  Let $G$  be a generic
subset of $\calP$  with $p\in G$.  Then  
${}^\ome\!N^*[G]\subseteq N^*[G]$.

\pr It is enough to show that for every $i<\del$
$${}^\ome\!N^*_i[G]\subseteq N^*_{i+1}[G]\ .$$
Let $i<\del$.  By Definition 5.1, $|N^*_i|=\del$,
$N^*_i\in N^*_{i+1}$,  $\del\subseteq
N^*_{i+1}$,  $\del <\kap_0$  and
$\del^\ome=\ome$.  Then $|^\ome\!N^*_i[G]|=\del$.
Since $|N^*_i[G]|=\del$  and $\del^\ome=\del$
in a generic extension as well.  But
$N^*_{i+1}[G]$  is an elementary submodel of
$V_\lam [G]$ for $\lam$  big enough.  So ${}^\ome\!N^*_i
[G]\subseteq N^*_{i+1}[G]$.\hfill$\square$

\sect{7.~~Extensions of Elementary Submodels in
the Projection Forcing}

Let $\langle D^*,\to\rangle$  the suborder of
$\langle\calP,\le\rangle$  defined in Section 4.
Denote by $\pi$  the corresponding projection.

\proclaim Lemma 7.1.  Let $\langle N,N^*\rangle$
be a good pair, $p\in\calP$  a $*$-generic for
$\langle N,N^*\rangle$ and $G\subseteq \calP$
generic subset with $p\in G$.  Then
\item{(1)} $N^*[G]\cap V\supseteq N^*[\pi\tagg (G)]\cap
V$.
\item{(2)} ${}^\ome\!N^*[\pi\tagg (G)]\subseteq 
N^*[\pi\tagg (G)]$.

\pr (1) holds since $G,\pi\in N^*[G]$.  For (2)
just repeat the argument of Lemma
6.8.\hfill$\square$

Let $\langle N^*_i\mid i<\del\rangle$  be a
sequence witnessing goodness of $N^*$.  By Lemma
6.4 it is possible to pick a $*$-generic $p$
over $N^*$ to be a supremum of
$\langle p_i\mid i<\del\rangle$  so that each $p_i$  is
$*$-generic over $N^*_i$  and $p_i\in N^*_{i+1}$. 

Let us denote by $\calP\rhookup p$  the set of
conditions
$\{q\in \calP|$  all the coordinates of $q$ are
already mentioned in $p\}$.

\proclaim Lemma 7.2. $N^*[G]\cap V$  and
$N^*[\pi\tagg G]\cap V$  depend only on $G\cap
\calP\rhookup p$.

\pr Let $\tau\in N^*$  be a name of an element
of $V$.  Then for some $i<\del$  $\tau\in
N^*_i$.  Since $p_i$  is $*$-generic for
$N^*_i$, there will be $n$,  $\ell(p_i)\le
n<\ome$  and a set measure one $A_{in}$  in
$p_i$  corresponding to $mc_n(p_i)$  such that
for every $\nu\in A_{in}$
$p_i^\cap\langle\nu\rangle$ decides the value of
$\tau$.  But then there will be $q\in G\cap
\calP\rhookup p_i$  $q\ge
p_i^\cap\langle\nu\rangle$  for some $\nu$.  
This $q$  will decide $\tau$.\hfill$\square$

\proclaim Lemma 7.2.1.  $N^*[G]\cap
V=N^*[\pi\tagg (G)]\cap V$.

\pr Let $\tau\in N^*$  be a name of element of
$V$.  As in 7.2, there are $i<\del$  and
$A_{in}$ such that $p_i\cap \langle\nu\rangle$
decides the value of $\tau$  for every $\nu\in
A_{in}$.  The point is that $\langle p_i^\cap
\langle \nu\rangle\mid\nu\in A_{in}\rangle$
remains a maximal antichain above $p_i$  also
after projection.  It follows from the
definition of $\pi$.  Now it is easy to replace
$\tau$  by a name $\tau'$  depending only on
$\langle p_i^\cap \langle\nu\rangle\mid\nu\in
A_i\rangle$  and it will be a name in the projection
forcing a well.  Now $\tau'[6]=\tau'[\pi\tagg
G]$.\hfill$\square$

\proclaim Lemma 7.3.  For every $n\ge \ell (p)$
$$\cup\supp\nolimits_n(p)=\cup
(N^*\cap\kap_n^{+n+2})=\cup (\supp\nolimits_n(p)\cap
N^*)\ .$$

\pr For every $\alp \in N^*$ $\cap\kap_n^{+n+2}$
there is $\bet\in N^*\cap\supp_n(p)$,  $\bet\ge
\alp$  since $p$  is $*$-generic for $N^*$.  On
the other hand,
$\supp_n(p)=\bigcup\limits_{i<\del}\supp_n(p_i)$
and for every $i<\del$  $p_i\in N^*$.  Hence
also $\cup\supp_n(p_i)\in N^*$.\hfill$\square$

\proclaim Lemma 7.4.  Let $n\ge\ell (p)$.  For
every $\alp\in N^*$ $\cap \kap_n^{+n+2}$  there
are $\bet, \gam\in N^*\cap$  $\supp_n(p)$ such
that $\bet$  is the $\alp$-th member of
$C_\gam$,  where $C_\gam$  is from the canonical box
sequence for $\kap^{+n+2}$.  

\pr Follows from elementary of $N^*$ and $*$-genericity
of $p$  for
$N^*$.\hfill$\square$

This lemma shows that we actually can recover
Prikry sequences even for coordinates $\alp$
which are occupied in $p$  by irrelevant
ordinals.  So $\calP\rhookup p$  is actually the
tree Prikry forcing where the measure including
$N^*\cap \kap_n^{+n+2}$  stands on the level $n$ in
the tree. 

Let $t\in D^*$.  Pick some good model $N^*$
with $t\in N^*$.  Let $p$  be $*$-generic over
$N^*$  and $t\le p$.  There is $p'\in D^*$  above
$p$.  Let $q'\in D^*$ $q'\lft p'$,  $q'\not=
p'$.  Pick a good model $M^*\supseteq N^*$  with
$q'\in M^*$. Let $r\ge q'$  be $*$-generic over $M^*$.
Pick $s,s\in D^*$ such that $s'\lft s$,  $s'\ge r$, 
$s\ge p'$.

\proclaim Lemma 7.5.  Let $H$ be a generic
subset of $\langle D^*,\to\rangle$  with $s\in
H$ (or equivalently $s'\in H$). Then
\item{(1)} both $N^*[H], M^*[H]$  are elementary
submodels of $V_\lam [H]$  
\item{(2)} if $G$, $G'$  are generic subsets
of $\langle\calP,\le\rangle$ s.t. $\pi\tagg
G=\pi\tagg G'=H$  and $s\in G$,  $s'\in G'$.
then
$$N^*[G]\cap V=N^*[H]\cap V$$
$$M^*[G']\cap V=M^*[H]\cap V\ .$$
\item{(3)} ${}^\ome\!M^*[H]\subseteq M^*[H]$ and
${}^\ome\!N^*[H]\subseteq N^*[H]$.

\pr (1) is trivial. (2),and (3) follow from
6.1.\hfill$\square$ 

Notice that since $q'\lft p'$  and $q'\not= p'$,
$M^*[H]$ and $N^*[G]$ will disagree about a
source of certain common Prikry sequence. 

The final step will be to show that this implies
a disagreement for indiscernibles as well.  Let's
fix $N^*$,  $M^*$,  $p,p',q',r,s,s'$  and $H$.
Notice that $\{ x\in \calP\mid x\ge s$  and
every coordinate of $x$  appears in $s\}\cap H$
is $\calP\rhookup s$  generic.  Denote it by
$H_s$  and define $H_{s'}$  in the same fashion.
Basically, we use $H$  to give us a Prikry
sequence for the maximal coordinates in $s$  (or
$s'$).  Use $H_{s'}$  to generate a
$\calP\rhookup r$-generic set $H_r$  in the same
fashion.  Let also $H_s$  generates a
$\calP\rhookup p$-generic set $H_p$.   

Combining Lemmas 7.2 and 7.5 we obtain the
following.

\proclaim Lemma 7.6.
\item{(a)} $N^*[H]\cap V=N^*[H_p]\cap V$
\item{(b)} $M^*[H]\cap V=M^*[H_r]\cap V$.
\item{}

Sets $H_p$ and $H_r$  will be used below to
produce winning strategies in games of the type
of [Git2].  We'll show that $H_p$  produces the
least winning strategy for $N^*$  and hence
coincides with the assignment functions for
indiscernibles of $N^*$.  The same will be true
about $H_r$  and $M^*$. 

Let us work in $V[H]$.  We denote by $\langle
c_n\mid n<\ome \rangle$  the Prikry sequence for
normal measures over $\kap_n$'s.  Assume that
$V=\calK(\calF)$  the core model with the maximal
sequence $\calF$  of extenders.  Let us recall
the statement of the Covering Lemma.

Suppose that there is no inner model with a
strong cardinal.  Let $\calK(\calF)$  be the
core model. 

\subheading{The Mitchell Covering Lemma}

Let $X<H_\lam$  (for some $\lam\ge\kap^+$) be
such that
\item{(a)} ${}^\ome\!X\subseteq X$
\item{(b)} $|X|<\kap$
\item{(c)} $X\cap\kap$  is cofinal in $\kap$.
\item{}

Let $\pi :N\simeq X$  be the transitive
collapse.  Denote $\pi^{-1}(\kap)$  by $\okap$,
$\pi^{-1}(\calF)$  by $\overline\calF$.

Then there are 
\item{(a)} $\overline\rho <\overline\kap$.
\item{(b)} a $\overline\calF\rhookup \kap$-mouse
$\om$ 
\item{(c)} functions on triples of ordinals
$\oC,C$  subsets of $\overline \kap$ and $\kap$
respectively
\item{(d)} a $\calF\rhookup\kap$-mouse $m$
\item{(e)} a map $\pi^*:\om\to m$
such that the following holds.
\item{(1)} $\om$  is an iterated ultrapower of
a $\overline\calF\rhookup\overline\rho$-mouse
and $\oC$  is the set of indiscernibles for
$\om$  corresponding to the iteration.
\item{(2)} $H^\om_\okap =N\cap
H^{\calK(\overline\calF)}_\okap$
\item{(3)} $\calP(\okap)\cap\calK(\overline\calF)\cap
N\subseteq \om$
\item{(4)} $\pi^*\supseteq\pi\rhookup (H_\okap)^{\calK
(\overline\calF)}$
\item{(5)} for every  $(\alp,\bet,\gam)\in\dom\oC$
$C(\pi^*(\alp),\pi^*(\bet),\pi^*(\gam))=\pi\tagg
(\oC(\alp, \bet,\gam))$ and $\dom C=\pi^*\,{}\!\tagg
\dom\oC$
\item{(6)} the canonical skolem function $\oh$
for $\om$  maps to the canonical skolem function
$h$ for $m$.
\item{(7)} $X\cap H_\kap\cap\calK(\calF)=h\tagg
(\pi\tagg(\orho);C)\cap H_\kap$
\item{(8)} $X\cap\calP(\kap)\cap\calK(\calF)\subseteq
h\tagg (\pi\tagg (\orho);C)$
\item{(9)} $H_\kap\cap \calK(\calF)\subseteq m$.
\item{}

Further we shall denote such $C$  and $h$ by
$C^X$  and $h^X$.  Also we shall confuse $C^X$
and $\cup \{ C^X(\alp,\bet,\gam)\mid \alp,\bet,\gam\in\dom
C^X\}$.  Elements of $C$  are called
indiscernibles.  

Notice that (1),(5),(6) imply that if $c\in rng
C^N$  then there is the unique triple $(\alp^N(c),
\bet^N(c),\gam^N(c))\in (h^N)\tagg c$  such that
$c\in C^N(\alp^N(c),\bet^N(c)$,  $\gam^N(c))$.

Let us now define a variant of the game
introduced in [Git2].

The game $\calG_\chi$ lasts $\le\chi$  moves,
where $\chi$  is an ordinal $<c_0$.  Actually,
it was used with $\chi=\ome_1+1$.

On stage $\del$ Player I picks a sequence 
$\oB_\del=\langle B_{\del n}\mid n<\ome\rangle$  such that
\item{(a)} For every $n<\ome$
$B_{\del_n}\in\calK(\calF)$,  it is a subset of
$\kap_n^{+n+2}$ of cardinality (in $\calK(\calF)$) $\le
\kap_n$.   
\item{(b)} There exists a function
$g\in\calK(\calF)$  such that for all (but
finitely many $n$'s) $B_{\del,n}\in g\tagg
(c_n)$.  

Notice that a particular case of (b) is
$\oB_\del=\langle B_{\del,n}\mid n<\ome
\rangle\in\calK(\calF)$.
\item{(c)} For every $n<\ome$
$B_{n,\del}\supseteq\bigcup_{\del'<\del}B_{n,\del'}$.

The answer of Player II is a sequence
$\otau^{\,\del}=\langle\otau^{\,\del}(n)\mid n<\ome
\rangle$ so that the conditions below hold
\item{(1)} For all but finitely many $n<\ome$
if $\langle\xi_i\mid i<\rho_{\del n}\rangle$ is
least in $\calK(\calF)$  enumeration of $B_{\del,n}$,
where $\rho_{\del,n}=\kap_n$  or $\rho_{\del,n}<\kap_n$
then $\otau^{\,\del}(n)=\langle\tau^\del
(n,i)\mid i<\rho'_{\del,n}\rangle$,  where
$$\rho'_{\del,n}=\cases{c_n(\alp_n)\ ,&if
$\rho_{\del,n}=\kap_n$\cr
\rho_{\del,n}\ ,&if
$\rho_{\del,n}<\kap_n$\ .\cr}$$
And for every $i,j<\rho'_n$
$$\tau^\del (n,i)<\tau^\del (n,j)\quad {\rm iff}\quad
\xi_i<\xi_j\ .$$
\item{(2)} For every $\del'<\del$ there is
$n(\del',\del )<\ome$  such that for every $n\ge
n(\del',\del)$  
$$\{\tau^\del (n,i)\mid i<\rho'_{\rho,n}\}\supseteq
\{\tau^{\del'}(n,j)\mid j<\rho'_{\del',n}\}$$
and
$\tau^\del (n,i)=\tau^{\del'}(n,j)$ implies that
the $i$-th element of $B_{\del,n}$  and the $j$-th
element of $B_{\del',n}$  are the same.

Player II loses if at some stage he has no legal
moves.

The idea behind this definition is as follows. 
The first player picks indexes of measures and
the second is supposed to find indiscernible (or
what he thinks are indiscernibles) for these
measures.  Then the first player increases the
set of indexes of measures and the second
chooses indiscernibles for the bigger set.  He
is supposed to respect his previous choices
almost everywhere.  If he succeeds in going on
long enough he wins.

Notice, that if Player I restricts himself to
moves that are in some fixed model $N$,  then II
has a winning strategy.  He just plays
indiscernibles of $N$.

Let's use $H_p$  to define a winning strategy
$\sig$  for II in the game restricted to $N^*[H]$.
The strategy will be positional, i.e. will
depend only on the last move.  So suppose that
stage $\del$  Player I plays a sequence $\oB_\del=\langle
B_{\del,n}\mid n<\ome\rangle\in N^*[H]$ such
that each $B_{\del n}\in\calK (\calF)$,  is a
subset of $\kap_n^{+n+2}$  of cardinality $\le\kap_n$
(in $\calK (\calF)$ or equivalently in $V[H]$)
and there is a function $g\in\calK(\calF)$ such
that for all but finitely many $n$'s $B_{\del
,n}\in g\tagg (c_n)$.  By elementarity of
$N^*[H]$, we can assume that $g\in N^*$.     

\proclaim Lemma 7.7.  $\langle B_{\del, n}\mid
n<\ome\rangle\in\calK(\calF)$.

\pr Let $g(\xi_n)=B_{\del,n}$  for some
$\xi_n<c_n$  and every $n<\ome$.  Using the
Prikry condition of $\calP\rhookup p$  and the
fact that $\langle c_n\mid n<\ome \rangle$  in a
Prikry sequence for normal measures, it is not
hard to see that $\langle\xi_n\mid n<\ome\rangle\in
N^*$.  But then also $\langle B_{\del,n}\mid
n<\ome\rangle\in N^*$.\hfill$\square$

For every $n<\ome$  there is $\del_n<\kap_n^{+n+2}$
which codes $B_{\del,n}$  in $\calK(\calF)$.
Thus take $\rho_n$  to be the index of
$B_{\del,n}$ in $\calK(\calF)$ -- least
enumeration of subsets of $\kap_n^{+n+2}$  of
cardinality $\le\kap_n$.  Then, clearly,
$\langle\rho_n\mid n<\ome\rangle\in N^*$.  We
like $\sig(\langle B_{\del,n}\mid n<\ome\rangle
)$  be the sequence obtained by decoding the
Prikry sequence for the measures corresponding
to $\langle\rho_n\mid n<\ome\rangle$.  It may be
the case that $H_p$  just does not have such
Prikry sequence.  But then Lemma 7.4 will be
used to produce one.   

\proclaim Lemma 7.8.  In $N^*[H]$  there is a
Prikry sequence $\langle \chi_n\mid n<\ome\rangle$
for $\langle\rho_n\mid n<\ome\rangle$.

\pr Let $n\ge\ell (p)$.  Since $\rho_n\in
N^*\cap\kap_n^{+n+2}$, by Lemma 7.4 there will
be $\bet_n,\gam_n\in N^*\cap\supp_n(p)$  such
that $\bet_n$  is the $\rho_n$-th member of
$C_{\gam_n}$.  Since $\bet_n,\gam_n\in\supp_n(p)$, 
there will be in $H_p$,  $\xi_n$,  $\rho_n$  the
elements of the one element Prikry sequence
corresponding to measures of $\bet_n,\gam_n$.
Then for some $\chi_n$, $\xi_n$  will be
$\chi_n$-th element of $C_{\rho_n}$.  Now,
$\langle\chi_n\mid n\ge \ell (p)\rangle$  will
be the desired Prikry sequence.\hfill$\square$  

For every $n,\ome >n\ge\ell (p)$  let $\otau^{\,\del}(n)$
be the subset of $c_n^{+n+2}$  with index
$\chi_n$ in the $\calK(\calF)$-least enumeration
of subsets of $c_n^{+n+2}$  of cardinality $\le
c_n$.  Finally we define $\sig(\langle B_{\del_n}\mid
n<\ome\rangle)=\langle\otau^{\,\del}(n)\mid
n<\ome\rangle$.

The choice of $\langle\chi_n\mid n<\ome\rangle$
insures that $\langle \otau_\del (n)\mid
n<\ome\rangle$  is a legal move in the game,
i.e. the conditions (1) and (2) of its
definition are satisfied.  This completes the
definition of $\sig$.

Let us call a strategy $\mu$  a least strategy
if for every other strategy $\mu'$, for every
first move of Player I consisting of some
ordinals $\langle \{\bet_n\}\mid n<\ome\rangle$
$\mu (<\{\bet_n\}\mid
n<\ome\rangle)(\ell)\le\mu'(\langle\{\bet_n\mid
n<\ome\rangle)(\ell)$  for all but finitely many
$\ell$'s.

\proclaim Lemma 7.9.  $\sig$  is the least
possible winning strategy for games $\calG_\chi$
restricted to $N^*[H]$ with $\chi >\ome_1$.

\pr Suppose otherwise.  Let $\sig'$  be another
strategy and suppose that for some sequence
of ordinals $\obet =\langle \bet_n\mid
n<\ome\rangle$,  $\bet_n<\kap_n^{+n+2}$  $(n<\ome)$

$\sig'(\langle\{\bet_n\}\mid n<\ome\rangle)\not\ge
\sig(\langle \{\bet_n\}\mid n<\ome\rangle)$.

Denote $\sig'(\langle\{\bet_n\}\mid n<\ome\rangle)$
by $\onu^{\,0}=\langle\nu^0(n)\mid
n<\ome\rangle$  and $\sig (\langle\{\bet_n\}\mid
n<\ome \rangle)$  by
$\otau^{\,\obet}=\langle\tau^{\,\obet}(n)\mid
n<\ome \rangle$.  Assume for simplicity that for
every $n<\ome$
$\tau^{\,\obet}(n)>\nu^0(n)$  and $\ell (p)=0$.

We like to find now for each $n<\ome$  a measure
to which $\nu^0(n)$  corresponds and then to
proceed as in [Git2] toward the contradiction.

Fix $n<\ome$.  Let $\buildrul{}^\sim\under\nu\in
N^*$  be a $\calP$-name of $\nu^0(n)$.  Since
$p$ is $*$-generic over $N^*$  and
$\nu^0(n)<\kap_n^{+n+2}$, adding the Prikry
sequence up to level $n$  in $p$ already decides
the value of
$\buildrul{}^\sim\under\nu$. $\buildrul{}^\sim\under\nu$
can be viewed as function of $n$ variables and
it will represent an ordinal $<\kap_n^{+n+2}$ in
the iterated ultrapower by $mc_0(p)$, $mc_1(p),\ldots
mc_n(p)$.  This ordinal gives the index of the
measure corresponding to $\nu_0(n)$.\hfill$\square$

Since the playing indiscernibles for $N^*[H]$ is
also the least winning strategy for Player II,
we obtain the following  

\proclaim Lemma 7.10.  $\sig$  is almost equal
to the indiscernible strategy, which means that
for every ordinal move $\langle \{\bet_n\}\mid
n<\ome\rangle$  of Player 1 for all but finitely
many $n$'s.

$\sig (\langle \{\bet_n\}\mid n<\ome\rangle)(k)$  is
the indiscernible for $\bet_k$  for all but finitely
many $k$'s.

Using the same arguments for $M^*[H]$, $r$,  and
the fact that $H_p$  and $H_r$  disagree about
the source of some Prikry sequence in
$N^*[H]\cap M^*[H]$  we obtain the following.

\proclaim Theorem 7.11.  In $V[H]$ there are
elementary submodels which disagree about common
$\ome$-sequence of indiscernibles. 

[Git2] implies then to deduce the following corollary.

\proclaim Corollary 7.12.  The game $\calG_\chi$
is undetermined in $V[H]$.

\vskip1truecm
\references{75}	

\ref{[Git1]} M. Gitik, The strength of the
failure of SCH, Ann. of Pure and Appl. Logic 51
(1991), 215-240.
\smallskip
\ref{[Git-Mit]} M. Gitik and W. Mitchell, Indiscernible
sequences for extenders and the Singular Cardinal
Hypothesis. 
\smallskip
\ref{[Git2]} M. Gitik, On measurable cardinals
violating GCH, to appear in Ann. of Pure and
Appl. Logic. 
\smallskip
\ref{[Sh]} S. Shelah, Cardinal Arithmetic, to appear.
\smallskip
\ref{[Git-Mag1]} M. Gitik and M. Magidor, The
Singular Cardinal Hypothesis Revisited, in MSRI
Conf., 1991, 243-279. 
\smallskip
\ref{[Git-Mag2]} M. Gitik and M. Magidor, Extender
Based Forcing Notions, to appear in JSL.

\end